\documentclass{article}
\usepackage{amssymb}
\usepackage{amsfonts}
\usepackage{amsmath}

\input{tcilatex}
\begin{document}

\title{Global variational solutions to a class of fractional spde's on
unbounded domains}
\author{Marco Dozzi$^{\ast }$, Rim Touibi$^{\ast }$, Pierre- A. Vuillermot$%
^{\ast ,\ast \ast }$ \\
UMR-CNRS 7502, Inst. \'{E}lie Cartan de Lorraine, Nancy, France$^{\ast }$,\\
Grupo de F\'{\i}sica Matem\'{a}tica, GFMUL, Faculdade de Ci\^{e}ncias,\\
Universidade de Lisboa, 1749-016 Lisboa, Portugal$^{\ast \ast }$}
\date{}
\maketitle

\begin{abstract}
In this article we prove new results regarding the existence and the
uniqueness of global variational solutions to Neumann initial-boundary value
problems for a class of non-autonomous stochastic parabolic partial
differential equations. The equations we consider are defined on unbounded
open domains in Euclidean space satisfying certain geometric conditions, and
are driven by a multiplicative noise derived from an infinite-dimensional
fractional Wiener process characterized by a sequence of Hurst parameters 
\textsf{H=}$\left( \mathsf{H}_{i}\right) _{i\in \mathbb{N}^{+}}\subset
\left( \frac{1}{2},1\right) $. These parameters are in fact subject to
further constraints that are intimately tied up with the nature of the
nonlinearity in the stochastic term of the equations, and with the choice of
the functional spaces in which the problem at hand is well-posed. Our method
of proof rests on compactness arguments in an essential way.
\end{abstract}

\section{Outline and statement of the main result}

Many applications ranging from the worlds of engineering and finance to the
natural sciences call for a mathematical modelling in terms of stochastic
differential equations. In particular, there have been many works devoted to
the analysis of the ultimate behavior of solutions to stochastic partial
differential equations of parabolic type which specifically occur in
population dynamics, population genetics, nerve pulse propagation and
related topics, to name only a few (see, e.g., \cite{vuillermot} for a brief
account of some of those works and the references therein). Moreover, there
have also been several more recent articles dealing with the analysis of
solutions to various types of semilinear parabolic stochastic partial
differential equations driven either by a Brownian noise, or by a fractional
noise with Hurst parameter $\mathsf{H}\in \left( \frac{1}{2},1\right) $
(see, e.g., \cite{alvarezlopezprivault}, \cite{balantudor}, \cite%
{dozzikolvovskalopez1}-\cite{dozzilopez}, and the plethora of references
therein, particularly \cite{lvduan}). While these works have been primarily
centered around questions of global existence, uniqueness and blowup in
finite time, there have also been investigations essentially motivated by
issues in financial mathematics devoted to the analysis of problems that
involve a mixture of a Brownian noise with a fractional noise,\textit{\ }%
within the realm of both ordinary and partial stochastic differential
equations (see, e.g., \cite{guerranualart}, \cite{mishurashevchenko1}-\cite%
{mishuraralshev} and the references therein).

It is our purpose here to contribute further to the analysis of some of the
above questions by proving the existence, and in some cases the uniqueness,
of global variational solutions to Neumann initial-boundary value problems
associated with a class of non-autonomous stochastic parabolic partial
differential equations defined on certain regions of Euclidean space. In
what follows we assume that all the functional spaces are real, use the
standard notations for the usual spaces of Lebesgue integrable functions and
their norms, and begin by defining the Wiener process that will generate the
fractional noise we need. Let $D\subset \mathbb{R}^{d}$ be an unbounded open
domain where $d\in \mathbb{N}^{+}$, and let us consider a linear,
self-adjoint, positive, non-degenerate trace-class operator $C$ in $L^{2}(D)$
whose eigenfunctions and eigenvalues we write $\left( e_{i}\right) _{i\in 
\mathbb{N}^{+}}$ and $\left( \lambda _{i}\right) _{i\in \mathbb{N}^{+}}$,
respectively. Let $\mathsf{H=}\left( \mathsf{H}_{i}\right) _{i\in \mathbb{N}%
^{+}}$ be a sequence of Hurst parameters satisfying $\mathsf{H}_{i}\in
\left( \frac{1}{2},1\right) $ for every $i$, and for $t\in \mathbb{R}%
_{0}^{+} $ let $\left( B^{\mathsf{H}_{i}}(t)\right) _{i\in \mathbb{N}^{+}}$
be a sequence of one-dimensional, independent fractional Brownian motions
defined on the complete probability space $\left( \Omega ,\mathcal{F},%
\mathbb{P}\right) $ and starting at the origin. Assuming that $\left(
e_{i}\right) _{i\in \mathbb{N}^{+}}$ constitutes an orthonormal basis of $%
L^{2}(D)$, we then define the $L^{2}(D)$-valued fractional Wiener process $%
W^{\mathsf{H}}\left( .,t\right) $ by%
\begin{equation}
W^{\mathsf{H}}\left( .,t\right) :=\sum_{i=1}^{+\infty }\sqrt{\lambda _{i}}%
e_{i}\left( .\right) B^{\mathsf{H}_{i}}(t).  \label{fractionalwiener}
\end{equation}%
This series converges strongly in $L^{2}(D)$ $\mathbb{P}$-a.s. by virtue of
the basic properties of $B^{\mathsf{H}_{i}}(t)$, the boundedness of the
sequence $\mathsf{H}$ and the fact that $C$ is trace-class. From this we
conclude that (\ref{fractionalwiener}) defines a centered Gaussian process
whose covariance is entirely determined by $C$, that is,%
\begin{eqnarray*}
&&\mathbb{E}\left( \left( W^{\mathsf{H}}\left( .,s\right) ,v\right)
_{2}\left( W^{\mathsf{H}}\left( .,t\right) ,\hat{v}\right) _{2}\right) \\
&=&\frac{1}{2}\sum_{i=1}^{+\infty }\left( s^{2\mathsf{H}_{i}}+t^{2\mathsf{H}%
_{i}}-\left\vert s-t\right\vert ^{2\mathsf{H}_{i}}\right) \left(
Cv,e_{i}\right) _{2}\left( \hat{v},e_{i}\right) _{2}
\end{eqnarray*}%
for all $s,t$ $\in \mathbb{R}_{0}^{+}$ and all $v,\hat{v}\in L^{2}(D)$,
where $\left( .,.\right) _{2}$ denotes the standard inner product in this
space and $\mathbb{E}$ the expectation functional on $\left( \Omega ,%
\mathcal{F},\mathbb{P}\right) $ (see, e.g., \cite{nuarasca} for a discussion
of various basic properties of fractional Brownian motion). For $T\in 
\mathbb{R}^{+\text{ }}$arbitrary, we then consider the class of stochastic
initial-boundary value problems formally given by%
\begin{eqnarray}
du(x,t) &=&(\limfunc{div}\left( k(x,t)\nabla _{x}u(x,t))+g(u(x,t)\right)
)dt+h(u(x,t))W^{\mathsf{H}}(x,dt),  \notag \\
(x,t) &\in &D\times \left( 0,T\right] ,  \notag \\
u(x,0) &=&\varphi (x),\text{ \ \ }x\in D,  \label{stochasticproblems} \\
\frac{\partial u(x,t)}{\partial n(k)} &=&0,\text{ \ \ }(x,t)\in \partial
D\times \left( 0,T\right]  \notag
\end{eqnarray}%
where $\partial D=\overline{D}\setminus D$ stands for the boundary of $D$
and the last line for the conormal derivative of $u$ relative to the
matrix-valued function $k$. Moreover, $g$ and $h$ are real-valued while $%
\varphi $ is a random initial condition. Regarding these functions we will
need the following hypotheses:

\bigskip

(K) The function $k:$ $D\times \left( 0,T\right] \mapsto \mathbb{R}^{d^{2}}$
is Lebesgue-mesurable and we have $k_{i,j}(.)=k_{j,i}(.)$ for all $i,j\in
\left\{ 1,...,d\right\} $. Moreover, there exist constants $\underline{k},%
\overline{k}\in \mathbb{R}^{+}$ such that the inequalities%
\begin{equation}
\underline{k}\left\vert y\right\vert ^{2}\leq \left( k(x,t)y,y\right) \leq 
\overline{k}\left\vert y\right\vert ^{2}  \label{bounds}
\end{equation}%
hold uniformly in $(x,t)\in $ $D\times \left( 0,T\right] $ for all $y\in 
\mathbb{R}^{d}$, where $\left( .,.\right) $ and $\left\vert .\right\vert $
denote the standard Euclidean inner product and the related norm in $\mathbb{%
R}^{d}$, respectively.

\bigskip

(L) The functions $g,h:\mathbb{R\mapsto R}$ are Lipschitz continuous.
Moreover, the derivative $h^{\prime }$ exists, is H\"{o}lder continuous with
exponent $\gamma \in \left( 0,1\right] $ and bounded. In addition we impose
that%
\begin{equation}
\text{\underline{H}}\in \left( \frac{1}{\gamma +1},1\right)
\label{constraint}
\end{equation}%
where \underline{H}\textsf{\ }$:=\inf_{i\in \mathbb{N}^{+}}\mathsf{H}_{i}$.

\bigskip

(I) The initial condition $\varphi $ is an $L^{2}(D)$-valued random variable.

\bigskip

Finally, whereas the above properties of the operator $C$ are sufficient to
define $W^{\mathsf{H}}$, they are not quite strong enough to allow us to
prove the result we are looking for. Recalling that $C$ is necessarily an
integral transform whose generating kernel we denote by $\kappa $, we still
impose the following spectral condition:

\bigskip

(C) We have 
\begin{equation}
x\mapsto \int_{D}dy\left\vert \kappa \left( x,y\right) \right\vert ^{2}\in
L^{\infty }(D)  \label{boundedness}
\end{equation}

and 
\begin{equation}
\sum_{i=1}^{+\infty }\sqrt{\lambda _{i}}\left\Vert e_{i}\right\Vert _{\infty
}<+\infty .  \label{spectralcondition}
\end{equation}

\bigskip

As to the consistency of this hypothesis, we simply remark that we have
indeed $e_{i}\in L^{\infty }(D)$ as a consequence of (\ref{boundedness}),
which follows easily from the eigenvalue equation $Ce_{i}=\lambda _{i}e_{i}$%
, Schwarz inequality and the fact that $\left\Vert e_{i}\right\Vert _{2}=1$.

We will also need some functional spaces in order to define the notion of
variational solution we are interested in. To this end we fix once and for
all an $\alpha \in \left( 1-\underline{\mathsf{H}},\frac{\gamma }{\gamma +1}%
\right) $ where \underline{H} satisfies (\ref{constraint}), and introduce
the Banach space $\mathcal{B}^{\alpha ,2}\left( \left[ 0,T\right]
;L^{2}(D)\right) $ of all Lebesgue measurable mappings $u:\left[ 0,T\right]
\mapsto L^{2}(D)$ endowed with the norm%
\begin{equation}
\left\Vert u\right\Vert _{\alpha ,2,T}^{2}:=\left( \sup_{t\in \left[ 0,T%
\right] }\left\Vert u(t)\right\Vert _{2}\right) ^{2}+\int_{0}^{T}dt\left(
\int_{0}^{t}d\tau \frac{\left\Vert u(t)-u(\tau )\right\Vert _{2}}{\left(
t-\tau \right) ^{\alpha +1}}\right) ^{2}<+\infty \text{.\label{firstnorm}}
\end{equation}%
Furthermore, let $H^{1}\left( D\times \left( 0,T\right) \right) $ be the
isotropic Sobolev space consisting of all real-valued functions $v\in
L^{2}\left( D\times \left( 0,T\right) \right) $ that possess distributional
derivatives $v_{x_{i}},v_{\tau }\in L^{2}\left( D\times \left( 0,T\right)
\right) $, endowed with the inner product%
\begin{eqnarray}
\left( v_{1},v_{2}\right) _{1,2,T} &:&=\int_{D\times \left( 0,T\right)
}dxd\tau v_{1}(x,\tau )v_{2}(x,\tau )  \notag \\
&&+\sum_{i=1}^{d}\int_{D\times \left( 0,T\right) }dxd\tau v_{1,x_{i}}(x,\tau
)v_{2,x_{i}}(x,\tau )  \notag \\
&&+\int_{D\times \left( 0,T\right) }dxd\tau v_{1_{,}\tau }(x,\tau
)v_{2_{,}\tau }(x,\tau )  \label{innerproduct}
\end{eqnarray}%
and the corresponding norm%
\begin{equation}
\left\Vert v\right\Vert _{1,2,T}=\left( v,v\right) _{1,2,T}^{\frac{1}{2}}%
\text{.}  \label{norm}
\end{equation}%
We note that for any function $v\in H^{1}\left( D\times \left( 0,T\right)
\right) $ which does not depend on time we have $v\in H^{1}\left( D\right) $%
, the usual Sobolev space on $D$ whose norm we denote by $\left\Vert
.\right\Vert _{1,2}$. Moreover, let $H^{1}\left( (0,T)\right) $ be the
Sobolev space of functions defined on the time interval $(0,T)$. From these
definitions it follows immediately that if $v\in H^{1}\left( D\right) $ and $%
\eta \in H^{1}\left( (0,T)\right) $ then $v\otimes \eta \in H^{1}\left(
D\times \left( 0,T\right) \right) $. This allows us to consider the vector
space of all finite linear combinations of such tensor products as an inner
product space with respect to (\ref{innerproduct}), and we write $\mathcal{H}%
\left( D\times \left( 0,T\right) \right) $ for its completion in $%
H^{1}\left( D\times \left( 0,T\right) \right) $ with respect to (\ref{norm}%
). About this Hilbert space we shall prove in Section 2 the existence of the
trace embedding%
\begin{equation}
\mathcal{H}\left( D\times \left( 0,T\right) \right) \rightarrow L^{2}\left(
D\times \left\{ \tau =t\right\} \right)  \label{tracebis}
\end{equation}%
valid for each $t\in \left[ 0,T\right] $ and every $d\in \mathbb{N}^{+}$,
which is important to formulate the following notion of variational solution:

\bigskip

\textbf{Definition 1.} We say the $H^{1}(D)$-valued random field $\left(
u_{V,\varphi }\left( .,t\right) \right) _{t\in \left[ 0,T\right] }$ defined
on $\left( \Omega ,\mathcal{F},\mathbb{P}\right) $ is a \textit{global
variational solution to (\ref{stochasticproblems})} if the following two
conditions hold:

(1) We have $u_{V,\varphi }\in L^{2}\left( 0,T;H^{1}(D)\right) \cap $ $%
\mathcal{B}^{\alpha ,2}\left( \left[ 0,T\right] ;L^{2}(D)\right) $ $\mathbb{P%
}$-a.s., which means that the relations%
\begin{equation}
\int_{0}^{T}d\tau \left\Vert u_{V,\varphi }\left( .,\tau \right) \right\Vert
_{1,2}^{2}=\int_{0}^{T}d\tau \left( \left\Vert u_{V,\varphi }\left( .,\tau
\right) \right\Vert _{2}^{2}+\left\Vert \nabla u_{V,\varphi }\left( .,\tau
\right) \right\Vert _{2}^{2}\right) <+\infty  \label{normbis}
\end{equation}%
and%
\begin{equation*}
\left\Vert u_{V,\varphi }\right\Vert _{\alpha ,2,T}<+\infty
\end{equation*}%
hold $\mathbb{P}$-a.s. In the sequel we shall sometimes write $\left\Vert
u_{V,\varphi }\right\Vert _{L^{2}\left( 0,T;H^{1}(D)\right) \text{ }}$for
norm (\ref{normbis}).

(2) The integral relation%
\begin{eqnarray}
\left( v(.,t),u_{V,\varphi }\left( .,t\right) \right) _{2} &=&\left(
v(.,0),\varphi \right) _{2}+\int_{0}^{t}d\tau \left( v_{\tau }(.,\tau
),u_{V,\varphi }\left( .,\tau \right) \right) _{2}  \notag \\
&&-\sum_{i,j=1}^{d}\int_{0}^{t}d\tau \left( v_{x_{i}}(.,\tau
),k_{i,j}(.,\tau )u_{V,\varphi ,x_{j}}\left( .,\tau \right) \right) _{2} 
\notag \\
&&+\int_{0}^{t}d\tau \left( v(.,\tau ),g(u_{V,\varphi }\left( .,\tau \right)
)\right) _{2}  \notag \\
&&+\int_{0}^{t}\int_{D}dxv(x,\tau )h(u_{V,\varphi }\left( x,\tau \right) )W^{%
\mathsf{H}}(x,d\tau )  \label{variasolution}
\end{eqnarray}%
holds $\mathbb{P}$-a.s. for each $v\in \mathcal{H}\left( D\times \left(
0,T\right) \right) $ and every $t\in \left[ 0,T\right] $, where $x\mapsto
v(x,t)\in L^{2}(D)$ stands for the trace of $v$ in the sense of (\ref%
{tracebis}), and where the stochastic integral with respect to (\ref%
{fractionalwiener}) is defined as%
\begin{eqnarray}
&&\int_{0}^{t}\int_{D}dxv(x,\tau )h(u_{V,\varphi }\left( x,\tau \right) )W^{%
\mathsf{H}}(x,d\tau )  \notag \\
&:&=\sum_{i=1}^{+\infty }\sqrt{\lambda _{i}}\int_{0}^{t}\left( v(.,\tau
),h(u_{V,\varphi }\left( .,\tau \right) )e_{i}\right) _{2}B^{\mathsf{H}%
_{i}}(d\tau ).  \label{stochasticintegral}
\end{eqnarray}

\bigskip

In order to make sense out of each term in (\ref{variasolution}) and prove
the existence of such a solution we will still need the following geometric
hypothesis on the domain $D$:

\bigskip

(D) There exists a sequence $(D_{N})_{N\in \mathbb{N}^{+}}$of open sets such
that for every $N$ we have $D_{N}\subset D_{N+1}\subset D$ and the
compactness of the embedding%
\begin{equation}
H^{1}\left( D_{N}\right) \hookrightarrow L^{2}\left( D_{N}\right) .
\label{compactembedding}
\end{equation}%
Moreover, let $\mathcal{B}\subset $ $H^{1}\left( D\right) $ be any bounded
set. Then for each $\epsilon >0$ there exists $N_{\epsilon }\in \mathbb{N}%
^{+}$ such that%
\begin{equation}
\int_{D\setminus D_{N_{\epsilon }}}dx\left\vert u(x)\right\vert ^{2}<\epsilon
\label{smallness}
\end{equation}%
for every $u\in \mathcal{B}$.

\bigskip

\textsc{Remark.} It is essential to observe here that any domain satisfying
(D) is necessarily of finite Lebesgue measure as a consequence of the theory
of tesselations developed in Chapter 6 of \cite{adams}, and more importantly
that (D) implies the compactness of the embedding $H^{1}\left( D\right)
\hookrightarrow L^{2}\left( D\right) $ as it amounts to a geometric
condition implying that $D$ becomes rapidly narrow at infinity, as the two
examples discussed at the end of this section will show. It is also easily
verified that each term in (\ref{variasolution}) is well defined and finite $%
\mathbb{P}$-a.s. as a consequence of all of the above hypotheses. In
particular, we may conclude that (\ref{stochasticintegral}) is an infinite
sum of one-dimensional, pathwise, generalized Stieltjes integrals which
defines a real-valued random variable as a consequence of Hypothesis (C) and
of the fact that $h$ is Lipschitz continuous. We shall dwell a bit more on
this and on related properties of the stochastic integral in the Appendix.
We note that Problem (\ref{stochasticproblems}) was thoroughly analyzed in 
\cite{nualartvuillerm} in case $D$ is a bounded domain satisfying the cone
condition and with a single Hurst parameter in (\ref{fractionalwiener}).

\bigskip

Under these conditions our main result is the following:

\bigskip

\textbf{Theorem. }\textit{Assume that Hypotheses (K), (L), (I), (C) and (D)
hold. Then Problem (\ref{stochasticproblems}) possesses a global variational
solution }$u_{V,\varphi }$\textit{. Moreover, if }$h$ \textit{is an affine
function then }$u_{V,\varphi }$\textit{\ is the unique solution to (\ref%
{stochasticproblems}).}

\bigskip

In order to prove this result we shall organize the remaining part of this
article in the following way: in Section 2 we first state the existence of $%
u_{V,\varphi }$ when the test functions in (\ref{variasolution}) are
independent of time, that is, with $v\in H^{1}\left( D\right) $. We then
extend the statement to the case of an approximating family of test
functions in $\mathbb{R}^{d+1}$ consisting of finite linear combinations of
the tensor products we alluded to above, and eventually to all test
functions $v\in \mathcal{H}\left( D\times \left( 0,T\right) \right) $ by a
suitable density argument and by invoking the trace embedding of the form (%
\ref{tracebis}) which we prove by elementary means. In Section 2 we also
illustrate our results with two examples, while we prove some crucial
estimates for the stochastic integral in the Appendix by means of a
necessary modification of the theory set forth in \cite{nuarasca} due to the
fact that we are dealing there with an infinite sequence of Hurst parameters
satisfying (\ref{constraint}). The method rests in an essential way on a
particular case of an inequality proved by Garsia, Rodemich and Rumsey in 
\cite{garsiaetal}, and on Minkowski's integral inequality (see, e.g.,
Appendix A in \cite{stein}).

\section{Proof of the main result and two examples}

We begin with the following result:

\bigskip

\textbf{Proposition 1.} \textit{Assume that Hypotheses} \textit{(K), (L),
(I), (C) and (D) hold. Then there exists }$u_{V,\varphi }\in L^{2}\left(
0,T;H^{1}(D)\right) \cap $ $\mathcal{B}^{\alpha ,2}\left( \left[ 0,T\right]
;L^{2}(D)\right) $ \textit{such that the integral relation (\ref%
{variasolution}) holds for every }$v\in H^{1}\left( D\right) $, \textit{that
is,}%
\begin{eqnarray}
\left( v,u_{V,\varphi }\left( .,t\right) \right) _{2} &=&\left( v,\varphi
\right) _{2}-\sum_{i,j=1}^{d}\int_{0}^{t}d\tau \left(
v_{x_{i}},k_{i,j}(.,\tau )u_{V,\varphi ,x_{j}}\left( .,\tau \right) \right)
_{2}  \notag \\
&&+\int_{0}^{t}d\tau \left( v,g(u_{V,\varphi }\left( .,\tau \right) )\right)
_{2}  \notag \\
&&+\int_{0}^{t}\int_{D}dxv(x)h(u_{V,\varphi }\left( x,\tau \right) )W^{%
\mathsf{H}}(x,d\tau )  \label{variasolutionter}
\end{eqnarray}%
$\mathbb{P}$\textit{-a.s. for every }$t\in \left[ 0,T\right] $.

\bigskip

\textbf{Proof. }Since the theory developed in Chapter 6 of \cite{adams}
implies the compactness of the embedding $H^{1}(D)\hookrightarrow L^{2}(D)$,
the result follows from a direct adaptation of all the arguments put forward
in the first part of \cite{nualartvuillerm}. \ \ $\blacksquare $

\bigskip

The next intermediary result is:

\bigskip

\textbf{Lemma 1. }\textit{Assume that Hypotheses} \textit{(K), (L), (I), (C)
and (D) hold and let }$u_{V,\varphi }$ \textit{be} \textit{the random field
of Proposition 1. Then for all finite linear combinations }$\hat{v}$\textit{%
\ of functions of the form }$v\otimes \eta $ \textit{where }$v\in
H^{1}\left( D\right) $\textit{\ and }$\eta \in H^{1}\left( \left( 0,T\right)
\right) $ \textit{we have }%
\begin{eqnarray}
\left( \hat{v}(.,t),u_{V,\varphi }\left( .,t\right) \right) _{2} &=&\left( 
\hat{v}(.,0),\varphi \right) _{2}+\int_{0}^{t}d\tau \left( \hat{v}_{\tau
}(.,\tau ),u_{V,\varphi }\left( .,\tau \right) \right) _{2}  \notag \\
&&-\sum_{i,j=1}^{d}\int_{0}^{t}d\tau \left( \hat{v}_{x_{i}}(.,\tau
),k_{i,j}(.,\tau )u_{V,\varphi ,x_{j}}\left( .,\tau \right) \right) _{2} 
\notag \\
&&+\int_{0}^{t}d\tau \left( \hat{v}(.,\tau ),g(u_{V,\varphi }\left( .,\tau
\right) )\right) _{2}  \notag \\
&&+\int_{0}^{t}\int_{D}dx\hat{v}(x,\tau )h(u_{V,\varphi }\left( x,\tau
\right) )W^{\mathsf{H}}(x,d\tau )  \label{approxrelation}
\end{eqnarray}%
$\mathbb{P}$\textit{-a.s. for every }$t\in \left[ 0,T\right] $.

\bigskip

\textbf{Proof. }By linearity it is sufficient to prove that (\ref%
{approxrelation}) holds for $\hat{v}=v\otimes \eta $. The easiest way out is
to start with the first integral on the right-hand side. Using (\ref%
{variasolutionter}) we obtain%
\begin{eqnarray}
&&\int_{0}^{t}d\tau \left( \hat{v}_{\tau }(.,\tau ),u_{V,\varphi }\left(
.,\tau \right) \right) _{2}  \notag \\
&=&\int_{0}^{t}d\tau \eta ^{\prime }(\tau )\left( v,u_{V,\varphi }\left(
.,\tau \right) \right) _{2}  \notag \\
&=&\left( \hat{v}(.,t),\varphi \right) _{2}-\left( \hat{v}(.,0),\varphi
\right) _{2}  \notag \\
&&-\sum_{i,j=1}^{d}\int_{0}^{t}d\tau \eta ^{\prime }(\tau )\int_{0}^{\tau
}d\sigma \left( v_{x_{i}},k_{i,j}(.,\sigma )u_{V,\varphi ,x_{j}}\left(
.,\sigma \right) \right) _{2}  \notag \\
&&+\int_{0}^{t}d\tau \eta ^{\prime }(\tau )\int_{0}^{\tau }d\sigma \left(
v,g(u_{V,\varphi }\left( .,\sigma \right) )\right) _{2}  \notag \\
&&+\int_{0}^{t}d\tau \eta ^{\prime }(\tau )\int_{0}^{\tau
}\int_{D}dxv(x)h(u_{V,\varphi }\left( x,\sigma \right) )W^{\mathsf{H}%
}(x,d\sigma ).  \label{substitution}
\end{eqnarray}%
We then integrate by parts the last three terms with respect to $\tau $ and
reintroduce $\hat{v}$ whenever possible to obtain respectively%
\begin{eqnarray}
&&\sum_{i,j=1}^{d}\int_{0}^{t}d\tau \eta ^{\prime }(\tau )\int_{0}^{\tau
}d\sigma \left( v_{x_{i}},k_{i,j}(.,\sigma )u_{V,\varphi ,x_{j}}\left(
.,\sigma \right) \right) _{2}  \notag \\
&=&\eta (t)\sum_{i,j=1}^{d}\int_{0}^{t}d\tau \left( v_{x_{i}},k_{i,j}(.,\tau
)u_{V,\varphi ,x_{j}}\left( .,\tau \right) \right) _{2}  \notag \\
&&-\sum_{i,j=1}^{d}\int_{0}^{t}d\tau \left( \hat{v}_{x_{i}}(.,\tau
),k_{i,j}(.,\tau )u_{V,\varphi ,x_{j}}\left( .,\tau \right) \right) _{2}
\label{intparts}
\end{eqnarray}%
and%
\begin{eqnarray}
&&\int_{0}^{t}d\tau \eta ^{\prime }(\tau )\int_{0}^{\tau }d\sigma \left(
v,g(u_{V,\varphi }\left( .,\sigma \right) )\right) _{2}  \notag \\
&=&\eta (t)\int_{0}^{t}d\tau \left( v,g(u_{V,\varphi }\left( .,\tau \right)
)\right) _{2}-\int_{0}^{t}d\tau \left( \hat{v}(.,\tau ),g(u_{V,\varphi
}\left( .,\tau \right) )\right) _{2}  \label{intpartsbis}
\end{eqnarray}%
for the deterministic integrals, while we get 
\begin{eqnarray}
&&\int_{0}^{t}d\tau \eta ^{\prime }(\tau )\int_{0}^{\tau
}\int_{D}dxv(x)h(u_{V,\varphi }\left( x,\sigma \right) )W^{\mathsf{H}%
}(x,d\sigma )  \notag \\
&=&\eta (t)\int_{0}^{t}\int_{D}dxv(x)h(u_{V,\varphi }\left( x,\tau \right)
)W^{\mathsf{H}}(x,d\tau )  \notag \\
&&-\int_{0}^{t}\int_{D}dx\hat{v}(x,\tau )h(u_{V,\varphi }\left( x,\tau
\right) )W^{\mathsf{H}}(x,d\tau )  \label{intpartster}
\end{eqnarray}%
for the stochastic integral. The substitution of (\ref{intparts})-(\ref%
{intpartster}) into (\ref{substitution}) then leads to the desired result,
after having lumped together the three terms containing the factor $\eta (t)$
and used there (\ref{variasolutionter}) once again. \ \ $\blacksquare $

\bigskip

We infer from the preceding considerations that for any $v\in \mathcal{H}%
\left( D\times \left( 0,T\right) \right) $ there exist functions $\hat{v}%
_{n} $ satisfying (\ref{approxrelation}) for every $n\in \mathbb{N}^{+}$
such that%
\begin{equation}
\left\Vert v-\hat{v}_{n}\right\Vert _{1,2,T}\rightarrow 0  \label{density}
\end{equation}%
as $n\rightarrow +\infty $, where $\left\Vert .\right\Vert _{1,2,T}$ is
given by (\ref{norm}). The first consequence of this is that we can already
approximate the three deterministic integrals in (\ref{variasolution}),
deferring to separate propositions the analysis of the remaining terms:

\bigskip

\textbf{Proposition 2.}\textit{\ Let }$v\in \mathcal{H}\left( D\times \left(
0,T\right) \right) $\textit{\ and let }$\left( \hat{v}_{n}\right) _{n\in 
\mathbb{N}^{+}}$\textit{\ be as in (\ref{density}). Then we have}%
\begin{equation}
\lim_{n\rightarrow +\infty }\int_{0}^{t}d\tau \left( \hat{v}_{n,\tau
}(.,\tau ),u_{V,\varphi }\left( .,\tau \right) \right)
_{2}=\int_{0}^{t}d\tau \left( v_{\tau }(.,\tau ),u_{V,\varphi }\left( .,\tau
\right) \right) _{2},  \label{limit1}
\end{equation}%
\begin{eqnarray}
&&\lim_{n\rightarrow +\infty }\sum_{i,j=1}^{d}\int_{0}^{t}d\tau \left( \hat{v%
}_{n,x_{i}}(.,\tau ),k_{i,j}(.,\tau )u_{V,\varphi ,x_{j}}\left( .,\tau
\right) \right) _{2}  \notag \\
&=&\sum_{i,j=1}^{d}\int_{0}^{t}d\tau \left( v_{x_{i}}(.,\tau
),k_{i,j}(.,\tau )u_{V,\varphi ,x_{j}}\left( .,\tau \right) \right) _{2}
\label{limit2}
\end{eqnarray}%
\textit{and}%
\begin{equation}
\lim_{n\rightarrow +\infty }\int_{0}^{t}d\tau \left( \hat{v}_{n}(.,\tau
),g(u_{V,\varphi }\left( .,\tau \right) )\right) _{2}=\int_{0}^{t}d\tau
\left( v(.,\tau ),g(u_{V,\varphi }\left( .,\tau \right) )\right) _{2}
\label{limit3}
\end{equation}%
$\mathbb{P}$\textit{-a.s. for every }$t\in \left[ 0,T\right] $.

\bigskip

\textbf{Proof.} Regarding (\ref{limit1}) we have%
\begin{eqnarray*}
&&\int_{0}^{t}d\tau \left\vert \left( v_{\tau }(.,\tau )-\hat{v}_{n,\tau
}(.,\tau ),u_{V,\varphi }\left( .,\tau \right) \right) _{2}\right\vert \\
&\leq &\left\Vert u_{V,\varphi }\right\Vert _{L^{2}\left(
0,T;H^{1}(D)\right) \text{ }}\left( \int_{D\times (0,T)}dxd\tau \left\vert
v_{\tau }(x,\tau )-\hat{v}_{n,\tau }(x,\tau )\right\vert ^{2}\right) ^{\frac{%
1}{2}} \\
&\leq &\left\Vert u_{V,\varphi }\right\Vert _{L^{2}\left(
0,T;H^{1}(D)\right) \text{ }}\left\Vert v-\hat{v}_{n}\right\Vert
_{1,2,T}\rightarrow 0
\end{eqnarray*}%
as $n\rightarrow +\infty $ $\mathbb{P}$-a.s., by using successively Schwarz
inequalities in $L^{2}(D)$ and on the time interval $(0,T)$ along with (\ref%
{normbis}) and (\ref{density}). The proof of (\ref{limit2}) follows from
similar arguments. Thus, noting that the matrix elements $k_{i,j}$ are
uniformly bounded for all $i,j\in \left\{ 1,...,d\right\} $ as a consequence
of (\ref{bounds}), we eventually get%
\begin{eqnarray*}
&&\sum_{i,j=1}^{d}\int_{0}^{t}d\tau \left\vert \left( v_{x_{i}}(.,\tau )-%
\hat{v}_{n,x_{i}}(.,\tau ),k_{i,j}(.,\tau )u_{V,\varphi ,x_{j}}\left( .,\tau
\right) \right) _{2}\right\vert \\
&\leq &c\left\Vert u_{V,\varphi }\right\Vert _{L^{2}\left(
0,T;H^{1}(D)\right) \text{ }}\sum_{j=1}^{d}\left( \int_{D\times
(0,T)}dxd\tau \left\vert v_{x_{j}}(x,\tau )-\hat{v}_{n,x_{j}}(x,\tau
)\right\vert ^{2}\right) ^{\frac{1}{2}} \\
&\leq &c\left\Vert u_{V,\varphi }\right\Vert _{L^{2}\left(
0,T;H^{1}(D)\right) \text{ }}\left\Vert v-\hat{v}_{n}\right\Vert
_{1,2,T}\rightarrow 0
\end{eqnarray*}%
as $n\rightarrow +\infty $ $\mathbb{P}$-a.s., where $c>0$ is a constant
depending only on $\underline{k}$, $\overline{k}$ and $d$. In a similar way
for (\ref{limit3}) we have%
\begin{eqnarray*}
&&\lim_{n\rightarrow +\infty }\int_{0}^{t}d\tau \left\vert \left( v(.,\tau )-%
\hat{v}_{n}(.,\tau ),g(u_{V,\varphi }\left( .,\tau \right) )\right)
_{2}\right\vert \\
&\leq &\left( c_{1}+c_{2}\left\Vert u_{V,\varphi }\right\Vert _{L^{2}\left(
0,T;H^{1}(D)\right) \text{ }}\right) \lim_{n\rightarrow +\infty }\left\Vert
v-\hat{v}_{n}\right\Vert _{1,2,T}=0
\end{eqnarray*}%
$\mathbb{P}$.-a.s. with $c_{1,2}>0$, which follows from the fact that $g$ is
Lipschitz continuous and $\left\vert D\right\vert <+\infty $ where $%
\left\vert D\right\vert $ stands for the Lebesgue measure of $D$. \ \ $%
\blacksquare $

\bigskip

In the sequel we will also need approximation properties for the remaining
terms in (\ref{variasolution}) that are similar to those of the above
proposition. However, the proof of those properties and eventually of the
main theorem will require two more preparatory results. The first one is:

\bigskip

\textbf{Lemma 2. }\textit{Let }$\mathcal{H}\left( D\times \left( 0,T\right)
\right) $\textit{\ be the Hilbert space defined in Section 1. Then there
exists the continuous trace embedding}%
\begin{equation}
\mathcal{H}\left( D\times \left( 0,T\right) \right) \rightarrow L^{2}\left(
D\times \left\{ \tau =t\right\} \right)  \label{embedding1}
\end{equation}%
\textit{valid for each }$t\in \left[ 0,T\right] $\textit{\ and every }$d\in 
\mathbb{N}^{+}$\textit{. Moreover, there also exists the continuous embedding%
}%
\begin{equation}
\mathcal{H}\left( D\times \left( 0,T\right) \right) \rightarrow \mathcal{B}%
^{\alpha ,2}\left( \left[ 0,T\right] ;L^{2}(D)\right)  \label{embedding2}
\end{equation}%
\textit{where }$\mathcal{B}^{\alpha ,2}\left( \left[ 0,T\right]
;L^{2}(D)\right) $\textit{\ is endowed with norm (\ref{firstnorm}).}

\bigskip

\textbf{Proof.} It is sufficient to prove the result for $H^{1}\left(
D\times \left( 0,T\right) \right) $. Let $v\in H^{1}\left( D\times \left(
0,T\right) \right) $ and let us write momentarily $\left\Vert .\right\Vert
_{1,2,H^{1}\left( (0,T)\right) }$ for the norm in $H^{1}\left( (0,T)\right) $%
. From (\ref{norm}) and Fubini's theorem we infer that $x$ $\mapsto
v(x,t)\in L^{2}\left( D\right) $ for almost every $t$, and more importantly
that $t\mapsto v(x,t)\in H^{1}\left( (0,T)\right) $ for almost every $x\in D$%
. Therefore, writing $\mathcal{C}\left( \left[ 0,T\right] \right) $ for the
space of all continuous functions defined on $\left[ 0,T\right] $ endowed
with the uniform norm $\left\Vert .\right\Vert _{\infty ,T}$ we have $%
t\mapsto v(x,t)\in \mathcal{C}\left( \left[ 0,T\right] \right) $ by virtue
of the embedding $H^{1}\left( (0,T)\right) \hookrightarrow \mathcal{C}\left( %
\left[ 0,T\right] \right) $. Consequently we obtain%
\begin{equation*}
\left\vert v(x,t)\right\vert \leq \left\Vert v(x,.)\right\Vert _{\infty
,T}\leq c\left\Vert v(x,.)\right\Vert _{1,2,H^{1}\left( (0,T)\right) }
\end{equation*}%
for each $t\in \left[ 0,T\right] $ and almost every $x\in D$ for some $c>0$,
and thereby the desired result%
\begin{equation*}
\int_{D}dx\left\vert v(x,t)\right\vert ^{2}\leq c^{2}\int_{D}dx\left\Vert
v(x,.)\right\Vert _{1,2,H^{1}\left( (0,T)\right) }^{2}\leq c^{2}\left\Vert
v\right\Vert _{1,2,T}^{2}
\end{equation*}%
according to (\ref{norm}). As for the proof of the second embedding, we note
that for every $v\in $ $H^{1}\left( D\times \left( 0,T\right) \right) $ and $%
t\geq \tau $ we have%
\begin{eqnarray*}
&&\left\Vert v(.,t)-v(.,\tau )\right\Vert _{2} \\
&\leq &\left( t-\tau \right) ^{\frac{1}{2}}\left( \int_{D\times
(0,T)}dxd\sigma \left\vert v_{\sigma }(x,\sigma )\right\vert ^{2}\right) ^{%
\frac{1}{2}}\leq \left( t-\tau \right) ^{\frac{1}{2}}\left\Vert v\right\Vert
_{1,2,T}
\end{eqnarray*}%
by (\ref{norm}), and therefore%
\begin{eqnarray*}
&&\left\Vert v\right\Vert _{\alpha ,2,T}^{2}=\left( \sup_{t\in \left[ 0,T%
\right] }\left\Vert v(.,t)\right\Vert _{2}\right) ^{2}+\int_{0}^{T}dt\left(
\int_{0}^{t}d\tau \frac{\left\Vert v(.,t)-v(.,\tau )\right\Vert _{2}}{\left(
t-\tau \right) ^{\alpha +1}}\right) ^{2} \\
&\leq &c\left\Vert v\right\Vert _{1,2,T}^{2}+\int_{0}^{T}dt\left(
\int_{0}^{t}d\tau \left( t-\tau \right) ^{-\alpha -\frac{1}{2}}\right)
^{2}\left\Vert v\right\Vert _{1,2,T}^{2}\leq c\left\Vert v\right\Vert
_{1,2,T}^{2}
\end{eqnarray*}%
by virtue of (\ref{embedding1}) and the fact that $\alpha <\frac{1}{2}$,
changing the value of $c$ whenever necessary. \ \ $\blacksquare $

\bigskip

The desired approximation property for the stochastic integral in (\ref%
{variasolution}) will now rest on (\ref{embedding2}) and on the following
estimate, which also shows that (\ref{stochasticintegral}) is H\"{o}lder
continuous with respect to the time variable:

\bigskip

\textbf{Lemma 3.} \textit{Let us consider the stochastic integral as defined
in (\ref{stochasticintegral})}. \textit{Then there exists a }$\mathbb{P}$-%
\textit{a.s.} \textit{finite, positive random variable }$r_{\alpha }^{%
\mathsf{H}}$ \textit{such that the estimate }%
\begin{eqnarray}
&&\left\vert \int_{t^{\ast }}^{t}\int_{D}dxv(x,\tau )h(u_{V,\varphi }\left(
x,\tau \right) )W^{\mathsf{H}}(x,d\tau )\right\vert  \notag \\
&\leq &r_{\alpha }^{\mathsf{H}}\left( 1+\left\Vert u_{V,\varphi }\right\Vert
_{\alpha ,2,T}\right) \left\Vert v\right\Vert _{\alpha ,2,T}\left\vert
t^{\ast }-t\right\vert ^{\frac{1}{2}}  \label{stochasticestimate}
\end{eqnarray}%
\textit{holds} $\mathbb{P}$\textit{-a.s. for every }$v\in \mathcal{B}%
^{\alpha ,2}\left( \left[ 0,T\right] ;L^{2}(D)\right) $\textit{\ and all }$%
t,t^{\ast }\in \left[ 0,T\right] $.

\bigskip

We defer the proof of this lemma to the Appendix, as it requires technical
tools regarding generalized Stieltjes integrals that do not directly pertain
to the main core of this article.

The expected approximation property of the stochastic integral is then the
following:

\bigskip

\textbf{Proposition 3. }\textit{Let }$v\in \mathcal{H}\left( D\times \left(
0,T\right) \right) $ \textit{and }$\left( \hat{v}_{n}\right) _{n\in \mathbb{N%
}^{+}}$\textit{\ be as in Proposition 2. Then we have}%
\begin{eqnarray*}
&&\lim_{n\rightarrow +\infty }\int_{0}^{t}\int_{D}dx\hat{v}_{n}(x,\tau
)h(u_{V,\varphi }\left( x,\tau \right) )W^{\mathsf{H}}(x,d\tau ) \\
&=&\int_{0}^{t}\int_{D}dxv(x,\tau )h(u_{V,\varphi }\left( x,\tau \right) )W^{%
\mathsf{H}}(x,d\tau )
\end{eqnarray*}%
$\mathbb{P}$\textit{-a.s. for every }$t\in \left[ 0,T\right] $.

\bigskip

\textbf{Proof.} From Relation (\ref{stochasticestimate}) and embedding (\ref%
{embedding2}) we have%
\begin{eqnarray*}
&&\left\vert \int_{0}^{t}\int_{D}dx\left( v(x,\tau )-\hat{v}_{n}(x,\tau
)\right) h(u_{V,\varphi }\left( x,\tau \right) )W^{\mathsf{H}}(x,d\tau
)\right\vert \\
&\leq &cr_{\alpha }^{\mathsf{H}}\left( 1+\left\Vert u_{V,\varphi
}\right\Vert _{\alpha ,2,T}\right) \left\Vert v-\hat{v}_{n}\right\Vert
_{1,2,T}
\end{eqnarray*}%
$\mathbb{P}$-a.s. for some $c>0$, hence the desired result from (\ref%
{density}). \ \ $\blacksquare $

\bigskip

The preceding considerations now lead to the following:

\bigskip

\textbf{Proof of the main result.} Let $v\in \mathcal{H}\left( D\times
\left( 0,T\right) \right) $ and let\textit{\ }$\left( \hat{v}_{n}\right)
_{n\in \mathbb{N}^{+}}$\textit{\ }be as in (\ref{density}). According to
Lemma 1 we have%
\begin{eqnarray}
\left( \hat{v}_{n}(.,t),u_{V,\varphi }\left( .,t\right) \right) _{2}
&=&\left( \hat{v}_{n}(.,0),\varphi \right) _{2}+\int_{0}^{t}d\tau \left( 
\hat{v}_{n,\tau }(.,\tau ),u_{V,\varphi }\left( .,\tau \right) \right) _{2} 
\notag \\
&&-\sum_{i,j=1}^{d}\int_{0}^{t}d\tau \left( \hat{v}_{n,x_{i}}(.,\tau
),k_{i,j}(.,\tau )u_{V,\varphi ,x_{j}}\left( .,\tau \right) \right) _{2} 
\notag \\
&&+\int_{0}^{t}d\tau \left( \hat{v}_{n}(.,\tau ),g(u_{V,\varphi }\left(
.,\tau \right) )\right) _{2}  \notag \\
&&+\int_{0}^{t}\int_{D}dx\hat{v}_{n}(x,\tau )h(u_{V,\varphi }\left( x,\tau
\right) )W^{\mathsf{H}}(x,d\tau )  \label{aproxrelationbis}
\end{eqnarray}%
$\mathbb{P}$-a.s. for each $t\in \left[ 0,T\right] $ and every $n$, and we
already know that each integral on the right-hand side of this expression
converges to the desired integral in (\ref{variasolution}) according to
Propositions 2 and 3. Let us now investigate the term on the left-hand side.
From Schwarz inequality in $L^{2}(D)$ and (\ref{firstnorm}) we have%
\begin{eqnarray*}
&&\left\vert \left( v(.,t)-\hat{v}_{n}(.,t),u_{V,\varphi }\left( .,t\right)
\right) _{2}\right\vert  \\
&\leq &\left\Vert u_{V,\varphi }\right\Vert _{\alpha ,2,T}\sup_{t\in \left[
0,T\right] }\left\Vert v(.,t)-\hat{v}_{n}(.,t)\right\Vert _{2} \\
&\leq &\left\Vert u_{V,\varphi }\right\Vert _{\alpha ,2,T}\left\Vert v-\hat{v%
}_{n}\right\Vert _{\alpha ,2,T}\leq c\left\Vert u_{V,\varphi }\right\Vert
_{\alpha ,2,T}\left\Vert v-\hat{v}_{n}\right\Vert _{1,2,T}\rightarrow 0
\end{eqnarray*}%
$\mathbb{P}$-a.s. as $n\rightarrow +\infty $ for some $c>0$ by virtue of (%
\ref{density}) and (\ref{embedding2}). Therefore we get%
\begin{equation*}
\lim_{n\rightarrow +\infty }\left( \hat{v}_{n}(.,t),u_{V,\varphi }\left(
.,t\right) \right) _{2}=\left( v(.,t),u_{V,\varphi }\left( .,t\right)
\right) _{2}
\end{equation*}%
$\mathbb{P}$-a.s. uniformly in $t\in \left[ 0,T\right] $, where $v(.,t)$
stands for the trace of $v$ in the sense of (\ref{embedding1}). Passing to
the limit in (\ref{aproxrelationbis}) we then obtain for the remaining term%
\begin{eqnarray}
\lim_{n\rightarrow +\infty }\left( \hat{v}_{n}(.,0),\varphi \right) _{2}
&=&\left( v(.,t),u_{V,\varphi }\left( .,t\right) \right)
_{2}-\int_{0}^{t}d\tau \left( v_{\tau }(.,\tau ),u_{V,\varphi }\left( .,\tau
\right) \right) _{2}  \notag \\
&&+\sum_{i,j=1}^{d}\int_{0}^{t}d\tau \left( v_{x_{i}}(.,\tau
),k_{i,j}(.,\tau )u_{V,\varphi ,x_{j}}\left( .,\tau \right) \right) _{2} 
\notag \\
&&-\int_{0}^{t}d\tau \left( v(.,\tau ),g(u_{V,\varphi }\left( .,\tau \right)
)\right) _{2}  \notag \\
&&-\int_{0}^{t}\int_{D}dxv(x,\tau )h(u_{V,\varphi }\left( x,\tau \right) )W^{%
\mathsf{H}}(x,d\tau )  \label{approxrelationter}
\end{eqnarray}%
$\mathbb{P}$-a.s. for every $t\in \left[ 0,T\right] $, hence%
\begin{equation*}
\lim_{n\rightarrow +\infty }\left( \hat{v}_{n}(.,0),\varphi \right)
_{2}=\left( v(.,0),u_{V,\varphi }\left( .,0\right) \right) _{2}
\end{equation*}%
by choosing $t=0$. But from (\ref{variasolutionter}) at $t=0$ we have $%
\left( v,u_{V,\varphi }\left( .,0\right) -\varphi \right) _{2}=0$ for every $%
v\in H^{1}(D)$ and thereby for every smooth and compactly supported $v\in $ $%
\mathcal{C}_{c}^{1}(D)$, the latter space being dense in $L^{2}(D)$.
Therefore $u_{V,\varphi }\left( .,0\right) -\varphi $ is orthogonal to $%
L^{2}(D)$, hence $u_{V,\varphi }\left( .,0\right) =\varphi $ so that%
\begin{equation*}
\lim_{n\rightarrow +\infty }\left( \hat{v}_{n}(.,0),\varphi \right)
_{2}=\left( v(.,0),\varphi \right) _{2}.
\end{equation*}%
The substitution of the preceding expression into (\ref{approxrelationter})
then proves (\ref{variasolution}) for every $v\in \mathcal{H}\left( D\times
\left( 0,T\right) \right) $. Finally, the proof that $u_{V,\varphi }$ is the
unique variational solution to (\ref{stochasticproblems}) satisfying (\ref%
{variasolutionter}) when $h$ is an affine function is identical to that
carried out in \cite{nualartvuillerm}. Therefore, let $\tilde{u}_{V,\varphi }
$ be another variational solution to (\ref{stochasticproblems}) satisfying (%
\ref{variasolution}). We then have $\tilde{u}_{V,\varphi }=u_{V,\varphi }$ $%
\mathbb{P}$-a.s. as a solution to (\ref{variasolutionter}), hence also as a
solution to (\ref{variasolution}). \ \ $\blacksquare $

\bigskip

\textsc{Remarks.} (1) The preceding considerations show that there are
actually two distinct types of equivalent variational solutions to (\ref%
{stochasticproblems}), to wit, one that satisfies (\ref{variasolutionter})
which involves test functions independent of time, and one satisfying (\ref%
{variasolution}). As was done in \cite{nualartvuillerm} where the class of
problems given by (\ref{stochasticproblems}) was defined and analysed on a
bounded domain $D$ satisfying the cone condition, we may refer to them as 
\textit{variational solutions of type I and II}, respectively. Moreover, for
the kind of equations considered in this article and to the best of our
knowledge, the problem of uniqueness of the solution in case of a general
nonlinearity $h$ in the noise term remains open.

(2) If $D$ is a bounded domain that satisfies the cone condition, the
natural space of test functions for variational solutions of type II is $%
H^{1}\left( D\times \left( 0,T\right) \right) $, as was shown in \cite%
{nualartvuillerm} by first using (\ref{variasolutionter}) to prove (\ref%
{variasolution}) for polynomial test functions, then for smooth test
functions by invoking the $\mathcal{C}^{1}$-version of Weierstrass'
approximation theorem, and eventually for functions in $H^{1}\left( D\times
\left( 0,T\right) \right) $ by a density argument. It is not possible to
argue in that way when $D$ is unbounded and merely satisfies Hypothesis (D).
Thus, it is still true that the space $\mathcal{C}^{1}\left( \overline{D}%
\times \left[ 0,T\right] \right) $ of all uniformly continuous and uniformly
bounded functions on $D\times \left( 0,T\right) $ is dense in $H^{1}\left(
D\times \left( 0,T\right) \right) $ if we assume in addition that $D$ and
therefore $D\times \left( 0,T\right) $ satisfy the so-called segment
condition (see, e.g. \cite{adams}). However, we can no longer approximate
functions in $\mathcal{C}^{1}\left( \overline{D}\times \left[ 0,T\right]
\right) $ by polynomials in a suitable topology to make such approximations
worthwile in the context of this article (an approximation of the functions
in $\mathcal{C}^{1}\left( \overline{D}\times \left[ 0,T\right] \right) $ by
polynomials is still possible if one endows $\mathcal{C}^{1}\left( D\times
\left( 0,T\right) \right) $ with the structure of a Fr\'{e}chet space, see,
e.g., Section 3 of Chapter 3 in \cite{kirillovgvi}, but this did not turn
out to be strong enough in our case). We therefore bypassed the difficulty
by constructing the space $\mathcal{H}\left( D\times \left( 0,T\right)
\right) $ in Section 1, but the determination of the relative size of this
space within $H^{1}\left( D\times \left( 0,T\right) \right) $ remains an
open problem at the present time.

\bigskip

We complete this section by illustrating the above results by means of two
examples. We begin with:

\bigskip

\textbf{Proposition 4.}\textit{\ Let us consider the two-dimensional domain}%
\begin{equation}
D=\left\{ \left( x_{1},x_{2}\right) \in \mathbb{R}^{2}:x_{1}>0,\text{ \ }%
0<x_{2}<b(x_{1})\right\}  \label{domain}
\end{equation}%
\textit{where the positive boundary curve }$b\in \mathcal{C}^{1}\left( 
\mathbb{R}_{0}^{+}\right) $ \textit{is decreasing, has a bounded derivative
and satisfies}%
\begin{equation}
\lim_{s\rightarrow +\infty }\frac{b\left( s+\epsilon \right) }{b\left(
s\right) }=0  \label{decrease}
\end{equation}%
\textit{for every }$\epsilon >0$\textit{. Then the conditions of Hypothesis
(D) in Section 1 hold.}

\bigskip

\textbf{Proof.} For every $N\in \mathbb{N}^{+}$ let us define the bounded
open set%
\begin{equation*}
D_{N}:=\left\{ \left( x_{1},x_{2}\right) \in \mathbb{R}^{2}:0<x_{1}<N\right%
\} \cap D\text{.}
\end{equation*}%
Then we have $D_{N}\subset D_{N+1}\subset D$ and $D_{N}$ satisfies the cone
condition, so that compact embedding (\ref{compactembedding}) holds by the
Rellich-Kondrachov theorem. Furthermore, the properties of the flow
associated with $D$ in the sense of Example 6.46 and Theorem 6.47 in \cite%
{adams} by means of the function $b$ satisfying (\ref{decrease}) imply that
for every $u\in H^{1}(D)$, the inequality%
\begin{equation*}
\int_{D\setminus D_{N}}dx\left\vert u(x)\right\vert ^{2}\leq c_{N}\left\Vert
u\right\Vert _{1,2}^{2}
\end{equation*}%
is valid with $\lim_{N\rightarrow +\infty }c_{N}=0$. Consequently, for each $%
\epsilon >0$ there exists $N_{\epsilon }\in \mathbb{N}^{+}$ such that (\ref%
{smallness}) holds for each $u$ belonging to a bounded set of $H^{1}(D)$. \
\ $\blacksquare $

\bigskip

We conclude that for a domain of the form (\ref{domain}) the embedding $%
H^{1}(D)\hookrightarrow L^{2}(D)$ is compact, so that our main theorem holds
true in this case. Typical examples of boundary functions satisfying the
above conditions include $b(x)=\exp \left[ -x^{1+\delta }\right] $ where $%
\delta >0$, whereas the mere exponential $b(x)=\exp \left[ -x\right] $ does
not belong to that class of curves. In fact, since $\ b$ is decreasing it
follows easily from (\ref{decrease}) that

\begin{equation*}
\lim_{s\rightarrow +\infty }\exp \left[ s\right] b(s)=0.
\end{equation*}%
Generally speaking, it is in fact the finiteness of the Lebesgue measure of $%
D$, along with the rapid decrease of the measure of the part of $D$ situated
outside the disk of radius $R$ centered at the origin as $R\rightarrow
+\infty $, which makes all this possible.

\bigskip

Our second example refers to the horn-shaped region discussed in \cite%
{adamsfournierbis}, for which our main theorem is also valid:

\bigskip

\textbf{Proposition 5.}\textit{\ Let us consider the three-dimensional domain%
}%
\begin{equation*}
D=\left\{ \left( x_{1},x_{2},x_{3}\right) \in \mathbb{R}^{3}:x_{3}>0,\text{
\ }0<x_{1}^{2}+x_{2}^{2}<b^{2}(x_{3})\right\}
\end{equation*}%
\textit{where the boundary curve }$b$\textit{\ satisfies the same hypotheses
as in Proposition 4 . Then the conditions of Hypothesis (D) in Section 1
hold.}

\bigskip

\textbf{Proof.} The argument is similar to that given in the proof of
Proposition 4 if we define%
\begin{equation*}
D_{N}:=\left\{ \left( x_{1},x_{2},x_{3}\right) \in \mathbb{R}%
^{3}:0<x_{3}<N\right\} \cap D
\end{equation*}%
for every $N\in \mathbb{N}^{+}$. Thus the compactness of embedding (\ref%
{compactembedding}) along with (\ref{smallness}) hold. \ \ $\blacksquare $

\bigskip

\textbf{Acknowledgements} The research of M. D. was supported in part by the
FCT of the Portuguese government under grant FCT UID/MAT/00208/2013. He
would like to thank Professors A. B. Cruzeiro, P.-A. Vuillermot and J. C.
Zambrini for their very kind invitation to the GFMUL, where this work was
completed. R. T. would like to extend her thanks and express her
appreciation and gratitude to the entire staff of the IECL for their warm
hospitality and strong support during the completion of her Ph.D.
dissertation there. The research of P.-A. V. was supported in part by the
FCT of the Portuguese government under grant PDTC/MAT-STA/0975/2014.

\bigskip

\textbf{Appendix: some remarks on generalized Stieltjes integrals and proof
of Lemma 3}

The following considerations constitute a necessary modification of the
theory developed in \cite{nuarasca}, needed to take into account the fact
that we are dealing here with an infinite sequence of Hurst parameters. For
every $i$ $\in \mathbb{N}^{+}$ let us introduce the random variable%
\begin{equation}
\Lambda _{\alpha }^{\mathsf{H}_{i}}=\frac{\sin \left( \pi \alpha \right) }{%
\pi }\sup_{0\leq t<t^{\ast }\leq T}\left\vert \frac{B^{\mathsf{H}_{i}}(t)-B^{%
\mathsf{H}_{i}}(t^{\ast })}{\left( t^{\ast }-t\right) ^{1-\alpha }}%
+(1-\alpha )\int_{t}^{t^{\ast }}d\tau \frac{B^{\mathsf{H}_{i}}(t)-B^{\mathsf{%
H}_{i}}(\tau )}{\left( \tau -t\right) ^{2-\alpha }}\right\vert
\label{randvarbis}
\end{equation}%
where $\alpha $ is the fixed parameter chosen in Section 1, for which we
have the inequality%
\begin{equation}
\Lambda _{\alpha }^{\mathsf{H}_{i}}\leq \sup_{0\leq t<t^{\ast }\leq T}\left( 
\frac{\left\vert B^{\mathsf{H}_{i}}(t^{\ast })-B^{\mathsf{H}%
_{i}}(t)\right\vert }{\left( t^{\ast }-t\right) ^{1-\alpha }}+(1-\alpha
)\int_{t}^{t^{\ast }}d\tau \frac{\left\vert B^{\mathsf{H}_{i}}(\tau )-B^{%
\mathsf{H}_{i}}(t)\right\vert }{\left( \tau -t\right) ^{2-\alpha }}\right) .
\label{randvarest}
\end{equation}%
For reasons that will soon be apparent, we first need to get appropriate
estimates for the moments of (\ref{randvarbis}) that are uniform in $i$ $\in 
\mathbb{N}^{+}$, and to achieve this suitable upper bounds for the
increments of $B^{\mathsf{H}_{i}}(t)$. Indeed the main result that will lead
to the proof of Lemma 3 is:

\bigskip

\textbf{Proposition A. 1.} \textit{Let }$\left( B^{\mathsf{H}_{i}}(t)\right)
_{i\in \mathbb{N}^{+}}$\textit{\ be the one-dimensional, independent
fractional Brownian motions introduced in Section 1, where the sequence }$%
\mathsf{H}=\left( \mathsf{H}_{i}\right) _{i\in \mathbb{N}^{+}}$\textit{\
satisfies (\ref{constraint}). Then for each }$i\in \mathbb{N}^{+}$\textit{\
we have}%
\begin{equation}
\sup_{i\in \mathbb{N}^{+}}\mathbb{E}\left\vert \Lambda _{\alpha }^{\mathsf{H}%
_{i}}\right\vert ^{p}<+\infty  \label{uniformity}
\end{equation}%
\textit{for every} $p\in \left[ 1,+\infty \right) $.

\bigskip

The proof of this proposition rests on the following particular case of the
Garsia-Rodemich-Rumsey inequality, obtained by rescaling the basic estimate
stated in Lemma 1.1 of \cite{garsiaetal} to establish its validity on $\left[
0,T\right] $ rather than just on $\left[ 0,1\right] $, and by applying the
particular choice of the functions involved made at the very beginning of 
\cite{hule}:

\bigskip

\textbf{Lemma A. 1. }\textit{Let }$f:\left[ 0,T\right] \mapsto \mathbb{R}$%
\textit{\ be continuous, let }$q\in \left[ 1,+\infty \right) $\textit{\ and }%
$\beta \in \left( \frac{1}{q},+\infty \right) $\textit{. If the integral
between parentheses in (\ref{garsiaetal}) is finite, then the inequality}%
\begin{equation}
\left\vert f(t^{\ast })-f(t)\right\vert \leq c_{\beta ,q,T}\left(
\int_{0}^{T}\int_{0}^{T}d\sigma d\tau \frac{\left\vert f(\sigma )-f(\tau
)\right\vert ^{q}}{\left\vert \sigma -\tau \right\vert ^{\beta q+1}}\right)
^{\frac{1}{q}}\left\vert t^{\ast }-t\right\vert ^{\beta -\frac{1}{q}}
\label{garsiaetal}
\end{equation}%
\textit{holds for all }$t,t^{\ast }\in \left[ 0,T\right] $\textit{, where}%
\begin{equation}
c_{\beta ,q,T}=c_{T}\left( \beta +\frac{1}{q}\right) \left( \beta -\frac{1}{q%
}\right) ^{-1}  \label{value}
\end{equation}%
\textit{and where }$c_{T}>0$\textit{\ depends only on }$T$\textit{.}

\bigskip

The upper bounds for the increments of $B^{\mathsf{H}_{i}}(t)$ we need turn
out to be provided by (\ref{garsiaetal}) for particular values of the
parameters. In all that follows we write $\epsilon $ for a fixed auxiliary
quantity which we can eventually express in terms of the constants $\alpha $%
, $\gamma $ and \b{H}:

\bigskip

\textbf{Lemma A. 2. }\textit{Let us fix }$\epsilon \in \left( 0,\frac{1}{%
\gamma +1}\right) $ \textit{where }$\gamma \in \left( 0,1\right] $\textit{\
is the constant of Hypothesis (L). Then for every }$i\in \mathbb{N}^{+}$ 
\textit{there exists a positive, }$\mathbb{P}$-\textit{a.s.} \textit{finite
random variable }$\Theta _{i,\epsilon ,T}$\textit{\ such that the inequality}%
\begin{equation}
\left\vert B^{\mathsf{H}_{i}}(t^{\ast })-B^{\mathsf{H}_{i}}(t)\right\vert
\leq c_{\gamma ,\epsilon ,T}\Theta _{i,\epsilon ,T}\left\vert t^{\ast
}-t\right\vert ^{\text{\underline{H}}-\epsilon }  \label{increstimate}
\end{equation}%
\textit{holds }$\mathbb{P}$-\textit{a.s.} \textit{for all }$t,t^{\ast }\in %
\left[ 0,T\right] $\textit{, where }$c_{\gamma ,\epsilon ,T}>0$\textit{\
depends only on }$\gamma $, $\epsilon $ \textit{and} $T$\textit{. Moreover
we have}%
\begin{equation}
\sup_{i\in \mathbb{N}^{+}}\mathbb{E}\left\vert \Theta _{i,\epsilon
,T}\right\vert ^{p}<+\infty \text{\label{uniformitybis}}
\end{equation}%
\textit{for every} $p\in \left[ 1,+\infty \right) $.

\bigskip

\textbf{Proof. }Let us first define the random variable%
\begin{equation*}
\tilde{\Theta}_{i,\epsilon ,T}:=\int_{0}^{T}\int_{0}^{T}d\sigma d\tau \frac{%
\left\vert B^{\mathsf{H}_{i}}(\sigma )-B^{\mathsf{H}_{i}}(\tau )\right\vert
^{\frac{2}{\epsilon }}}{\left\vert \sigma -\tau \right\vert ^{\frac{2\mathsf{%
H}_{i}}{\epsilon }}}
\end{equation*}%
and prove that%
\begin{equation}
\sup_{i\in \mathbb{N}^{+}}\mathbb{E}\left\vert \tilde{\Theta}_{i,\epsilon
,T}\right\vert ^{p}<+\infty  \label{uniformityquarto}
\end{equation}%
for every $p\in \left[ 1,+\infty \right) $. For any $r\in \left[ \frac{2}{%
\epsilon },+\infty \right) $ we have $p:=\frac{r\epsilon }{2}\in \left[
1,+\infty \right) $ so that on the one hand we obtain%
\begin{eqnarray}
\mathbb{E}\left\vert \tilde{\Theta}_{i,\epsilon ,T}\right\vert ^{p} &=&%
\mathbb{E}\left( \int_{0}^{T}\int_{0}^{T}d\sigma d\tau \frac{\left\vert B^{%
\mathsf{H}_{i}}(\sigma )-B^{\mathsf{H}_{i}}(\tau )\right\vert ^{\frac{2}{%
\epsilon }}}{\left\vert \sigma -\tau \right\vert ^{\frac{2\mathsf{H}_{i}}{%
\epsilon }}}\right) ^{\frac{r\epsilon }{2}}  \notag \\
&\leq &\left( \int_{0}^{T}\int_{0}^{T}d\sigma d\tau \left( \mathbb{E}\frac{%
\left\vert B^{\mathsf{H}_{i}}(\sigma )-B^{\mathsf{H}_{i}}(\tau )\right\vert
^{r}}{\left\vert \sigma -\tau \right\vert ^{r\mathsf{H}_{i}}}\right) ^{\frac{%
2}{r\epsilon }}\right) ^{\frac{r\epsilon }{2}}  \label{minkowski}
\end{eqnarray}%
from Minkowski's integral inequality. On the other hand, from the basic
properties of the one-dimensional fractional Brownian motion we have%
\begin{equation*}
\mathbb{E}\left\vert B^{\mathsf{H}_{i}}(\sigma )-B^{\mathsf{H}_{i}}(\tau
)\right\vert ^{r}\leq c_{r}\left\vert \sigma -\tau \right\vert ^{r\mathsf{H}%
_{i}}
\end{equation*}%
for some constant $c_{r}$ depending only on $r$, so that the substitution of
this expression into (\ref{minkowski}) annihilates the dependence in $%
\mathsf{H}_{i}$, thus leading to%
\begin{equation*}
\mathbb{E}\left\vert \tilde{\Theta}_{i,\epsilon ,T}\right\vert ^{p}\leq
c_{r}T^{2p}<+\infty
\end{equation*}%
uniformly in $i$, which is (\ref{uniformityquarto}). Let us now choose $q=%
\frac{2}{\epsilon }$ and $\beta _{i}=\mathsf{H}_{i}-\frac{\epsilon }{2}$ for
every $i\in \mathbb{N}^{+}$ in Lemma A. 1. Then clearly $q\in \left[
1,+\infty \right) $\textit{\ }and\textit{\ }$\beta _{i}\in \left( \frac{%
\epsilon }{2},+\infty \right) $ by virtue of (\ref{constraint}), and since (%
\ref{uniformityquarto}) implies in particular that $\tilde{\Theta}%
_{i,\epsilon ,T}<+\infty $ $\mathbb{P}$-a.s. we may thus apply (\ref%
{garsiaetal}) to obtain%
\begin{eqnarray*}
&&\left\vert B^{\mathsf{H}_{i}}(t^{\ast })-B^{\mathsf{H}_{i}}(t)\right\vert
\\
&\leq &c_{\beta _{i},q,T}\left( \int_{0}^{T}\int_{0}^{T}d\sigma d\tau \frac{%
\left\vert B^{\mathsf{H}_{i}}(\sigma )-B^{\mathsf{H}_{i}}(\tau )\right\vert
^{\frac{2}{\epsilon }}}{\left\vert \sigma -\tau \right\vert ^{\frac{2\mathsf{%
H}_{i}}{\epsilon }}}\right) ^{\frac{\epsilon }{2}}\left\vert t^{\ast
}-t\right\vert ^{\mathsf{H}_{i}-\epsilon } \\
&\leq &\tilde{c}_{\beta _{i},q,T}\left( \int_{0}^{T}\int_{0}^{T}d\sigma
d\tau \frac{\left\vert B^{\mathsf{H}_{i}}(\sigma )-B^{\mathsf{H}_{i}}(\tau
)\right\vert ^{\frac{2}{\epsilon }}}{\left\vert \sigma -\tau \right\vert ^{%
\frac{2\mathsf{H}_{i}}{\epsilon }}}\right) ^{\frac{\epsilon }{2}}\left\vert
t^{\ast }-t\right\vert ^{\text{\underline{H}}-\epsilon } \\
&=&\tilde{c}_{\beta _{i},q,T}\tilde{\Theta}_{i,\epsilon ,T}^{\frac{\epsilon 
}{2}}\left\vert t^{\ast }-t\right\vert ^{\text{\underline{H}}-\epsilon }
\end{eqnarray*}%
$\mathbb{P}$-a.s\textit{. }for all $t,t^{\ast }\in \left[ 0,T\right] $,
where $\tilde{c}_{\beta _{i},q,T}$ differs from $c_{\beta _{i},q,T}$ by a
trivial factor depending only on $T$ and \underline{H}. In order to get (\ref%
{increstimate}) it is thus sufficient to take%
\begin{equation}
\Theta _{i,\epsilon ,T}:=\tilde{\Theta}_{i,\epsilon ,T}^{\frac{\epsilon }{2}}
\label{randomvariable}
\end{equation}%
and prove that 
\begin{equation}
\sup_{i\in \mathbb{N}^{+}}\tilde{c}_{\beta _{i},q,T}<+\infty
\label{uniformityter}
\end{equation}%
along with (\ref{uniformitybis}). Ignoring the trivial dependence in $T$ and 
\underline{H} in $\tilde{c}_{\beta _{i},q,T}$ we first infer from (\ref%
{value}) that the simple estimate%
\begin{equation*}
\left( \beta _{i}+\frac{\epsilon }{2}\right) \left( \beta _{i}-\frac{%
\epsilon }{2}\right) ^{-1}=\mathsf{H}_{i}\left( \mathsf{H}_{i}-\epsilon
\right) ^{-1}\leq \left( \frac{1}{\gamma +1}-\epsilon \right) ^{-1}<+\infty
\end{equation*}%
holds uniformly in $i$ thanks to (\ref{constraint}) and our choice of $%
\epsilon $, which gives (\ref{uniformityter}). Finally, let us partition the
probability space as%
\begin{equation*}
\Omega =\left\{ \omega \in \Omega :\tilde{\Theta}_{i,\epsilon ,T}(\omega
)\leq 1\right\} \cup \left\{ \omega \in \Omega :\tilde{\Theta}_{i,\epsilon
,T}(\omega )>1\right\}
\end{equation*}%
and split the expectation functional accordingly. From (\ref{randomvariable}%
) it is then plain that%
\begin{equation*}
\sup_{i\in \mathbb{N}^{+}}\mathbb{E}\left\vert \Theta _{i,\epsilon
,T}\right\vert ^{p}=\sup_{i\in \mathbb{N}^{+}}\mathbb{E}\left\vert \tilde{%
\Theta}_{i,\epsilon ,T}\right\vert ^{\frac{p\epsilon }{2}}\leq 1+\sup_{i\in 
\mathbb{N}^{+}}\mathbb{E}\left\vert \tilde{\Theta}_{i,\epsilon
,T}\right\vert ^{p}<+\infty
\end{equation*}%
according to (\ref{uniformityquarto}) since $\frac{\epsilon }{2}\leq 1$,
which is the desired result. \ \ $\blacksquare $

\bigskip

The proof of (\ref{uniformity}) then follows from (\ref{uniformitybis}) and
from the substitution of (\ref{increstimate}) into (\ref{randvarest})
provided we impose a further restriction on the parameter $\epsilon $ to
make the singularities integrable:

\bigskip

\textbf{Proof of Proposition A. 1.} We first notice that 
\begin{equation*}
\text{\underline{H}}-1+\alpha <\frac{1}{\gamma +1}
\end{equation*}%
as a consequence of the conditions we imposed on these parameters in Section
1. Then, fixing $\epsilon \in \left( 0,\text{\underline{H}}-1+\alpha \right) 
$ we may substitute (\ref{increstimate}) into (\ref{randvarest}) to obtain%
\begin{equation*}
\Lambda _{\alpha }^{\mathsf{H}_{i}}\leq c_{\gamma ,\epsilon ,T}T^{\text{%
\underline{H}}-\epsilon -1+\alpha }\frac{\text{\underline{H}}-\epsilon
+\alpha }{\text{\underline{H}}-\epsilon -1+\alpha }\Theta _{i,\epsilon
,T}<+\infty
\end{equation*}%
$\mathbb{P}$-a.s. after an explicit integration, where the prefactor is
uniform in $i\in \mathbb{N}^{+}$. Therefore, (\ref{uniformity}) indeed
follows from (\ref{uniformitybis}). \ \ $\blacksquare $

\bigskip

\textsc{Remark.} To follow up on our remark preceding the statement of Lemma
A. 2., we see \textit{a posteriori} that we could have chosen for instance $%
\epsilon =\frac{1}{2}\left( \text{\underline{H}}-1+\alpha \right) $
throughout.

\bigskip

We are now ready for the following:

\bigskip

\textbf{Proof of Lemma 3.} We begin by estimating\textbf{\ }each integral on
the right-hand side of (\ref{stochasticintegral}). For every $i\in \mathbb{N}%
^{+}$ let us set%
\begin{equation*}
f_{i}(\tau ):=\left( v(.,\tau ),h(u_{V,\varphi }\left( .,\tau \right)
)e_{i}\right) _{2}.
\end{equation*}%
Then, from the basic estimate (4.11) in Proposition 4.1 of \cite{nuarasca}
regarding generalized Stieltjes integrals we infer that the inequality%
\begin{equation}
\left\vert \int_{t^{\ast }}^{t}f_{i}(\tau )B^{\mathsf{H}_{i}}(d\tau
)\right\vert \leq \Lambda _{\alpha }^{\mathsf{H}_{i}}\int_{t^{\ast
}}^{t}d\tau \left( \frac{\left\vert f_{i}(\tau )\right\vert }{\left( \tau
-t^{\ast }\right) ^{\alpha }}+\alpha \int_{t^{\ast }}^{\tau }d\rho \frac{%
\left\vert f_{i}(\tau )-f_{i}(\rho )\right\vert }{\left( \tau -\rho \right)
^{\alpha +1}}\right)  \label{stochasticestimatebis}
\end{equation}%
holds $\mathbb{P}$\textit{-}a.s., where $\Lambda _{\alpha }^{\mathsf{H}_{i}}$
is given by (\ref{randvarbis}). Furthermore, from Schwarz inequality and the
basic hypotheses of Section 1 we have%
\begin{eqnarray}
\left\vert f_{i}(\tau )\right\vert &\leq &c\left\Vert e_{i}\right\Vert
_{\infty }\left( 1+\left\Vert u_{V,\varphi }\left( .,\tau \right)
\right\Vert _{2}\right) \left\Vert v\left( .,\tau \right) \right\Vert _{2} 
\notag \\
&\leq &c\left\Vert e_{i}\right\Vert _{\infty }\left( 1+\left\Vert
u_{V,\varphi }\right\Vert _{\alpha ,2,T}\right) \left\Vert v\right\Vert
_{\alpha ,2,T}  \label{estimatequinto}
\end{eqnarray}%
where the second inequality follows from (\ref{firstnorm}), and similarly%
\begin{eqnarray}
&&\left\vert f_{i}(\tau )-f_{i}(\rho )\right\vert  \notag \\
&\leq &c\left\Vert e_{i}\right\Vert _{\infty }\left\Vert u_{V,\varphi
}\left( .,\tau \right) -u_{V,\varphi }\left( .,\rho \right) \right\Vert
_{2}\left\Vert v\right\Vert _{\alpha ,2,T}  \notag \\
&&+c\left\Vert e_{i}\right\Vert _{\infty }\left( 1+\left\Vert u_{V,\varphi
}\right\Vert _{\alpha ,2,T}\right) \left\Vert v\left( .,\tau \right)
-v\left( .,\rho \right) \right\Vert _{2}  \label{estimatesexto}
\end{eqnarray}%
for all $\rho ,\tau \in \left[ 0,T\right] $, for some constant $c>0$.
Consequently, on the one hand the substitution of (\ref{estimatequinto})
into the first integral on the right-hand side of (\ref%
{stochasticestimatebis}) gives%
\begin{equation}
\int_{t^{\ast }}^{t}d\tau \frac{\left\vert f_{i}(\tau )\right\vert }{\left(
\tau -t^{\ast }\right) ^{\alpha }}\leq c\left\Vert e_{i}\right\Vert _{\infty
}\left( 1+\left\Vert u_{V,\varphi }\right\Vert _{\alpha ,2,T}\right)
\left\Vert v\right\Vert _{\alpha ,2,T}\left\vert t^{\ast }-t\right\vert
^{1-\alpha }  \label{estimateseptimo}
\end{equation}%
by direct integration. On the other hand, the substitution of (\ref%
{estimatesexto}) into the second integral on the right-hand side of (\ref%
{stochasticestimatebis}), Schwarz inequality relative to the measure $d\tau $
and (\ref{firstnorm}) lead to%
\begin{equation}
\int_{t^{\ast }}^{t}d\tau \int_{t^{\ast }}^{\tau }d\rho \frac{\left\vert
f_{i}(\tau )-f_{i}(\rho )\right\vert }{\left( \tau -\rho \right) ^{\alpha +1}%
}\leq c\left\Vert e_{i}\right\Vert _{\infty }\left( 1+\left\Vert
u_{V,\varphi }\right\Vert _{\alpha ,2,T}\right) \left\Vert v\right\Vert
_{\alpha ,2,T}\left\vert t^{\ast }-t\right\vert ^{\frac{1}{2}}.
\label{estimateoitavo}
\end{equation}%
Therefore, with (\ref{estimateseptimo}) and (\ref{estimateoitavo}) into (\ref%
{stochasticestimatebis}) we obtain%
\begin{equation}
\left\vert \int_{t^{\ast }}^{t}f_{i}(\tau )B^{\mathsf{H}_{i}}(d\tau
)\right\vert \leq c\left\Vert e_{i}\right\Vert _{\infty }\Lambda _{\alpha }^{%
\mathsf{H}_{i}}\left( 1+\left\Vert u_{V,\varphi }\right\Vert _{\alpha
,2,T}\right) \left\Vert v\right\Vert _{\alpha ,2,T}\left\vert t^{\ast
}-t\right\vert ^{\frac{1}{2}}  \label{estimatenono}
\end{equation}%
for every $i\in \mathbb{N}^{+}$ since $1-\alpha >\frac{1}{2}$ according to
our original choice of this parameter. In order to prove (\ref%
{stochasticestimate}) it is therefore sufficient to show that%
\begin{equation}
\sum_{i=1}^{+\infty }\sqrt{\lambda _{i}}\left\Vert e_{i}\right\Vert _{\infty
}\Lambda _{\alpha }^{\mathsf{H}_{i}}<+\infty  \label{convergence}
\end{equation}%
$\mathbb{P}$\textit{-}a.s., for then%
\begin{eqnarray*}
&&\left\vert \int_{t^{\ast }}^{t}\int_{D}dxv(x,\tau )h(u_{V,\varphi }\left(
x,\tau \right) )W^{\mathsf{H}}(x,d\tau )\right\vert \\
&\leq &\left( c\sum_{i=1}^{+\infty }\sqrt{\lambda _{i}}\left\Vert
e_{i}\right\Vert _{\infty }\Lambda _{\alpha }^{\mathsf{H}_{i}}\right) \left(
1+\left\Vert u_{V,\varphi }\right\Vert _{\alpha ,2,T}\right) \left\Vert
v\right\Vert _{\alpha ,2,T}\left\vert t^{\ast }-t\right\vert ^{\frac{1}{2}}
\end{eqnarray*}%
according to (\ref{stochasticintegral}) and (\ref{estimatenono}), with the
obvious choice for $r_{\alpha }^{\mathsf{H}}$. But (\ref{convergence})
follows from Relation (\ref{uniformity}) of Proposition A.1. with $p=1$
since spectral condition (\ref{spectralcondition}) holds. \ \ $\blacksquare $

\bigskip

\bigskip

\bigskip

\bigskip

\bigskip

\end{document}